\documentclass[11pt]{amsart}

\usepackage{amsmath}
\usepackage{amssymb}
\usepackage{graphicx}
\newtheorem{theorem}{Theorem}[section]

\newtheorem{lemma}[theorem]{Lemma}
\newtheorem{proposition}[theorem]{Proposition}
\theoremstyle{definition}

\newcommand{\Proof}[1][]{{\par\smallskip\noindent{\bf Proof#1.\enspace }}}
\newcommand{\cqfd}{\hfill\rule{0.35em}{0.35em}\par\medskip\noindent}

\newcommand{\R}{\mathbb{R}}
\renewcommand{\div}{\operatorname{div}}

\def \R {I\!\!R}

\newcommand{\ba}{\begin{eqnarray}}
\newcommand{\ea}{\end{eqnarray}}

\newcommand{\basn}{\begin{eqnarray*}}
\newcommand{\easn}{\end{eqnarray*}}

\usepackage{xcolor}

\numberwithin{equation}{section}
\title[A Nonlinear problem with a weight ]{ A Nonlinear problem with a weight 
\\and a nonvanishing  boundary datum}

\author[Rejeb Hadiji]{Rejeb Hadiji}

\address[Rejeb Hadiji]{Universit\'e Paris-Est,
LAMA, Laboratoire d'Analyse et de Math\'ematiques Appliqu\'ees, \\
UMR 8050 du CNRS,\\
 61, avenue du G\'en\'eral de Gaulle, F-94010 Cr\'eteil -- France.}
\email{{\tt rejeb.hadiji@u-pec.fr}}



\begin{document}

\begin{abstract}
We consider the problem: 
$$\inf_{{u}\in {H}^{1}_{g}(\Omega),\|u\|_{q}=1} \int_{\Omega}{p(x)}\vert\nabla{u(x)}\vert^{2}dx-\lambda\int_{\Omega}\vert u(x)\vert^{2}dx$$
where $\Omega$ is a bounded domain in $\R^{n}$, ${n}\geq{4}$, $ p : \bar{\Omega}\longrightarrow \R$ is a given positive weight such that
$p\in H^{1}(\Omega)\cap C(\bar{\Omega})$, $0< c_1 \leq p(x) \leq c_2$,  $\lambda$ is a real
constant and $q=\frac{2n}{n-2}$ and $g$ a given positive boundary data. 
The goal of this present paper is to show that  minimizers do exist. We distinguish two cases, the first is solved by a convex argument while the second is not so straightforward and will be treated using 
the behavior of the weight near its minimum and the fact that the boundary datum is not zero.
\end{abstract}

\maketitle

\section{Introduction }

Let $\Omega$ be a  bounded domain in $\R^n$ of  class $C^{1}$, $n\geq 3$. Let us consider the minimization problem
\begin{eqnarray}\label{Pinf}
S_{0}(p,g)   = \inf_{{u}\in
{H}^{1}_{g}(\Omega),\|u\|_{q}=1}\int_{\Omega}{p(x)}\vert\nabla{u(x)}\vert^{2}dx  
\end{eqnarray}

\noindent where    
$${H}^{1}_{g}(\Omega) = \{ u\in {H}^{1} (\Omega) \,\hbox{s.t.}\, u = g \,\,\hbox{on}\,\, \partial \Omega\},$$
 $g\in H^{1\over 2} (\partial \Omega) \cap C(\partial \Omega )$ is a given boundary datum and $q=\frac{2n}{n-2}$ is the critical Sobolev exponent. 
 
 \noindent Note that it  is well known that  $H^{1}(\Omega)\hookrightarrow L^{r}(\Omega)$ is  continuous for any  $1 \leq r \leq  \frac{2n}{n-2}$. Moreover this embedding is compact for $1\leq r <  \frac{2n}{n-2}$.\\
We suppose that the weight  $p : \overline {\Omega} \rightarrow \,\R$ is a smooth function such that $0<c_1 \leq p(x)\leq c_2$ $\forall\,x\,\in\bar{\Omega}$ and 
$p$ is  in $H^1(\Omega)\cap C(\overline\Omega)$.

In this paper, we ask the question whenever the problem (\ref{Pinf}) has a minimizer. 
Note that if the infimum  (\ref{Pinf}) is achieved by some $u$ then we have
\ba \left\{
\begin{array}{lll} -\textsl{$\div$}(p(x)\nabla u)={ \Lambda u}^{q-{1}} &\textrm{in $\Omega$,}\\
\hspace{23mm}u > 0 &\textrm{in $\Omega$,}\\
\hspace{23mm}u=g  &\textrm{on $\partial \Omega$,}
\end{array}
\right. \label{Pequation1} \ea 
where $ \Lambda \in \R$ is  the Lagrange multiplier associated to the problem (\ref{Pinf}).

These  kind of  problems, which are  known to bear features of noncompactness are studied  by many authors.
First existence results for the  problem with a linear perturbation are due to Brezis-Nirenberg. Set 
 \begin{eqnarray}\label{Pinflambda}
S_{\lambda}(p,g)   = \inf_{{u}\in
{H}^{1}_{g}(\Omega),\|u\|_{q}=1}\int_{\Omega}{p(x)}\vert\nabla{u(x)}\vert^{2}dx  -\lambda\int_{\Omega}\vert u(x)\vert^{2}dx
\end{eqnarray}
  
They showed  that  if $g = 0$ and $p=1$, then $S_\lambda (1,0)$ is attainted as soon as  $S_\lambda (1,0) < S$ and this is the case if   $n\geq 4$, $0< \lambda < \lambda_1$, or $n =3$ and    $0 < \lambda^* < \lambda < \lambda_1$ where $\lambda_1$ is the first eigenvalue of $-\Delta$  and $\lambda^*$ depends on the domain, (see  \cite{6}). They showed also that   if $g\not\equiv 0$, $\lambda = 0$ and $p=1$ then  the infimum in  (\ref{Pinf}) is achieved, (see \cite{7}). Our approach uses their method.

In the case of $p=1$ and $g = 0$, Coron, Bahri and Coron exploited the topology of the domain. They  proved that equation 
\ba \left\{
\begin{array}{lll} - \Delta u ={ u}^{q-{1}} &\textrm{in $\Omega$,}\\
\hskip 6mm u > 0 &\textrm{in $\Omega$,}\\
\hskip 6mm u=0  &\textrm{on $\partial \Omega$,}
\end{array}
\right.\nonumber \ea 

\noindent  has a solution provided that the domain has nontrivial topology, (see \cite{8} and \cite{3}). 

\noindent We refer to \cite{13},  \cite{14} for the study of existence and multiplicity solutions of problem  (\ref{Pinf}) with the presence of a  smooth and positive  weight and  with homogeneous Dirichlet boundary condition. Nevertheless, in \cite{12}, it is shown that if $p$ is  discontinuous then  a solution of $S_0 (p,0)$ still exists.

\noindent In \cite{10}, the authors studied the minimization problem  on compact manifolds  in the case  $\lambda = 0$  with many variants.

\noindent For more general weights, depending on $x$ and on $u$, in a recent paper written  with Vigneron, we  showed  that in the case of homogeneous Dirichlet  boundary condition and in the presence of a linear perturbation the corresponding minimizing problem possesses a solution. The model of the  weight is $p(x,u) = \alpha + \vert x\vert^\beta \vert u\vert^k$ with positive parameters $\alpha$, $\beta$ and $k$. Note that in this case natural scalings appear and the answer depends on the ratio ${\beta \over k}$. For more details, we refer to \cite{2} and \cite{15}.

To motivate our problem, we briefly recall that it  is inspired by the study of the classical  Yamabe problem which has been the source of a large literature, 
(see for example  \cite{1}, \cite{3},  \cite{6}, \cite{8}, \cite{10} and \cite{16}), we refer to \cite{15} and the references therein for many recent developments in quasi-linear elliptic equations.

In this paper, we will assume that if $g\not\equiv 0$ having a constant sign and the weight  $p$ has a global  minimum $a\in \Omega$ such that  satisfies: 
\begin{eqnarray} \label{condition p}
p(x) \leq p_0 + \gamma\vert x-a\vert^{\alpha} \quad \forall x\in B(a,R)\subset \Omega, 
\end{eqnarray}
for constants  $\alpha > 1$, $\gamma > 0 $ and $R> 0$. 

The following  auxiliary linear  Dirichlet problem  will play an important role in  this paper:
 \begin{equation}\label{vvv}
\left\{
\begin{array}{llllll}
-\div(p\nabla v)=0 &\textrm{\quad in $\Omega$},\\
\hskip 17mm v=g  &\textrm{\quad on $\partial\Omega$}.
\end{array}
\right.
\end{equation}

 \subsection{Statement of the main result}

Our main result is  the following: 

\begin{theorem}\label{th1}
Let us assume that  the dimension $n\geq  3$  and $g\in H^{1\over 2} (\partial \Omega) \cap C(\partial \Omega )$ is a given boundary datum. Let $v$ be  the unique solution of (\ref{vvv}).  We have
\begin{enumerate}
\item 
Let $\vert\vert v\vert\vert_q < 1$ and let assume that $g\not\equiv 0$ and having a constant sign. Assume that $p$ has a global minimum $a\in \Omega$  that   satisfies (\ref{condition p}).
Then  for every $  n \in [3, 2\alpha + 2[$  the infimum  $S_{0}(p,g)$ is  achieved in $H^1_g(\Omega)$.
\item If $\vert \vert v\vert\vert_q \geq 1$ then for every $n\geq 3$ the infimum $S_{0}(p,g)$ is  achieved in $H^1_g(\Omega)$.
\end{enumerate}
\end{theorem}

The next proposition tell us that one has $\Sigma_g = \{ u \in {H}^{1}_{g}(\Omega),\|u\|_{q} =1\} \not= \emptyset$ which  ensures that $S_{0}(p,g)$ is well defined: 
\begin{proposition}
 Let  $g \in H^{1\over 2} (\partial \Omega) \cap C(\partial \Omega )$ be given boundary datum  and $v $ be the unique solution of (\ref{vvv}), we have 
 \begin{itemize}
\item If $\vert\vert v\vert\vert_q < 1$,then there is a bijection between $\Sigma_0$ and $\Sigma_g$. 
\item If $\vert\vert v\vert\vert_q \geq 1$, then $\Sigma_g \not= \emptyset$.
\end{itemize}
\label{proposition1}
\end{proposition}

Our problem depends on $\vert\vert v\vert\vert_q $. More precisely, we will use a convex argument to show that if $\vert\vert v\vert\vert_q \geq 1$ then the infimum (\ref{Pinf}) is achieved, 
while the case where  $\vert\vert v\vert\vert_q < 1$ is not so straightforward and  will be treated using the behavior of $p$ near its minimum  and the fact that $g$ has a constant sign. 
We will argue by contradiction, supposing that  minimizing sequence  converges weakly to some limit $u$. The fact that the boundary datum is not $0$ will give us  that $u$ is not identically $0$. 
Then, by using a suitable test functions, we will show equality (\ref{lemma1-2}) below which is  due to term of order $0$. After precise computations, we get strict inequality in (\ref{contradiction1})  which is due  to the next term in the same  expansion, which is lead to a contradiction.

Since  the nonlinearity of the problem is as stronger as $n$  is low,  it is rather surprising that the infimum is achieved for  lower dimensions $n\in [3,2\alpha + 2[$. 
Note that the presence of $p$ is  more significative if $\alpha >0$ is low. The compromise  is that  $n\in [3,2\alpha + 2[$. Remark  that if  $\alpha = 0$ then infimum of $p = p_0 + \gamma $ is not $p_0$.

\noindent For general boundary data $g$, we do not have  control over the normal derivative of a  solution of (\ref{Pequation1}) on the boundary of $\Omega$ and then, standard Pohozaev identity cannot be used.

\subsection{Structure of the paper} The paper is structured as follows: In  section 2 we give the notations and some preliminary results. 

\noindent In the next section, we state two results related to our main result namely, Theorem \ref{th2} which gives the sign of the Lagrange-multiplier associated to minimizers of $S_{0}(p, g)$ given by Theorem \ref{th1} and Theorem \ref{th3} which generalizes our main result in case of the presence of a linear perturbation.

\noindent In section 4, we will focus on the proof of Theorem \ref{th1}, which is the main result of this paper, it will be proved by a contradiction argument that spans the whole of this section. 

\noindent In section 5, we give the proof of Theorem \ref{th2}.

\noindent The last section is dedicated to  the problem  of  existence of minimizer in the presence of a linear perturbation and the proof of Theorem \ref{th3}.

\bigskip

\section{ Notations and  preliminary results }

Sobolev inequality says that there exist $M>0$ such that
\begin{equation*}
\int_{\Omega}p(x)|\nabla \phi|^{2}dx\geq M\left(\int_{\Omega}|\phi|^{q}dx\right)^{\frac{2}{q}}\quad\textrm{for all $\phi\,\in H_{0}^{1}(\Omega).$}
\end{equation*}

\noindent The best constant is defined by
$$
S_0 (p,0)=\inf_{\begin{array}{ll}u\in H_{0}^{1}(\Omega), ||u||_{q}=1\end{array}}\int_{\Omega}p(x)|\nabla u|^{2}dx.
$$

\noindent Set 
$$
S = S_0 (1,0) = \inf_{\begin{array}{ll}u\in H_{0}^{1}(\Omega), ||u||_{q}=1\end{array}}\int_{\Omega}|\nabla u|^{2}dx.
$$

\noindent We know that when the domain is $\R^n$, the constant $S_0 (1,0)$ is achieved   by the functions:
$$
U_{x_{0},\,\varepsilon}(x)=\left(\frac{\varepsilon}{\varepsilon^{2}+|x-x_{0}|^{2}}\right)^{\frac{n-2}{2}},\,\,\,x\,\in\,\R^{n}
$$
\noindent where $x_{0}\,\in\,\R^{n}$ and  $\varepsilon>0$, (see \cite{1}, \cite{6}, \cite{16}). 
Let us denote by  \begin{equation} \label{psi}
u_{x_{0},\varepsilon}(x)=U_{x_{0},\,\varepsilon}(x)\psi(x)
\end{equation}
where $\psi\,\in\,C^{\infty}(\R^{n})$, $\psi\equiv 1$ in $B(x_{0}, r)$ \,$\psi\equiv 0$ on $B(x_{0},\,2r)\subset \Omega$, $r>0$. We have
\begin{equation}\label{K1}
\int_{\Omega} p(x) |\nabla u_{x_{0},\varepsilon}|^{2}dx = p(x_{0}) K_1+ O(\varepsilon^{n-2}),
\end{equation}

\begin{equation}\label{K2}
\int_{\Omega}|u_{x_{0},\varepsilon}|^{q}dx= K_2 +O(\varepsilon^{n}),
\end{equation}

\noindent where $K_{1}$ and $K_{2}$ are positive constants with $\frac{K_{1}}{K_{2}^{2\over q}}=S$.

\noindent  We have also

\begin{equation*}
u_{x_{0},\,\varepsilon}\rightharpoonup 0 \quad\hbox{in}\quad  H_{0}^{1}(\Omega).
\end{equation*}
$$
-\Delta  U_{x_{0},\,\varepsilon}=c_{n}\,U_{x_{0},\,\varepsilon}^{q-1}\quad\hbox {in } \quad \R^{n}.
$$
\noindent It is well known that $S$ in never achieved for bounded domain, (see \cite{6}). 

\noindent In the the presence of the weight $p$ we have 
\begin{proposition}
\noindent Suppose that $a\in \Omega$ be a global minimum of $p$. 
Set $p_0  = p(a)$. If  $g = 0$, we have 
$S_{0} (p, 0)$ is never achieved and 
\begin{eqnarray*}
S_{0} (p, 0)=p_{0}\, S_{0} (1,0)=p_{0}\,S.
\end{eqnarray*}
\end{proposition}

\Proof
\noindent When $g=0$, the functions  $\frac{u_{a,\varepsilon}}{\|u_{a,\varepsilon}\|_{q}}$ are admissible test functions for $S_{0} (p, 0)$ and  we have as $\varepsilon\rightarrow 0$
\begin{eqnarray*}
p_0S \leq S_{0} (p, 0)&\leq& \int_{\Omega}p(x)\left|\nabla\frac{u_{a,\varepsilon}}{\|u_{a,\varepsilon}\|_{q}}\right|^{2}dx\\[\medskipamount] 
&=&p_{0}\,S+\int_{\Omega}(p(x)-p_{0})\left|\nabla\frac{u_{a,\varepsilon}}{\|u_{a,\varepsilon}\|_{q}}\right|^{2}dx  + o(1)\\[\medskipamount]
&=&p_{0}\,S +o(1).
\end{eqnarray*}
\noindent Passing to the limit $\varepsilon\rightarrow 0$  state that $S_{0}(p,\,0) = p_{0}S.$

\noindent This implies that $S_{0}(p,0)$ is not achieved. Indeed,  let us suppose that $S_{0}(p,0)$ is achieved by some $u$. Using the fact that  $S$ is never achived in bounded domains, we obtain
$$p_{0}S<p_{0}\int_{\Omega}|\nabla u|^{2}dx \leq \int_{\Omega}p(x)|\nabla u|^{2}dx=p_{0}S. $$
\noindent This leads to a contradiction.
 \cqfd
\subsection {The auxiliary Dirichlet problem}

The  linear  Dirichlet problem (\ref{vvv}) has a unique solution which solves  the following problem
\begin{eqnarray}\label{Pv}
 \min_{{v}\in{H}^{1}_{g}(\Omega)}\int_{\Omega}{p(x)}\vert\nabla{v(x)}\vert^{2}dx. 
\end{eqnarray}

\noindent Let us give now the proof of Proposition \ref{proposition1}: Recall that 
$$\Sigma_g = \{ u \in {H}^{1}_{g}(\Omega),\|u\|_{q} =1\}$$.
In the first case we can construct a bijection between $\Sigma_0$ and $\Sigma_g$. Indeed, 
let us  define, for $t$ in $\R$ and $u\in \Sigma_0$ the function
\begin{equation}\label{f(t)}
f(t) = \int_\Omega \vert tu + v\vert^q
\end{equation} 
since $f$ is smooth,  $f''(t)  = q(q-1) \int_\Omega \vert tu + v\vert^{q-2}u^2$,  $f(0) <1$ and $\lim_{t \rightarrow \infty} f(t) = \infty$, using the intermediate value theorem and the convexity of $f$ , we obtain,  for every  $u$ in $\Sigma_0$,  the  existence of  a unique $t(u)>0$ such  that $\vert\vert t(u)u +  v \vert\vert_q = 1$.\par

\noindent Let us denote by  $\varphi :  \Sigma_0\rightarrow \Sigma_g$  the function defined by $\varphi (u) = t(u)u + v$. 
Let $u_1$ and $u_2$ in $\Sigma_0$ such that $\varphi (u_1) = t(u_1)u_1 + v = \varphi (u_2) = t(u_2)u_2 + v$ , we have nesseceraly  $\vert\vert  t(u_1)u_1\vert\vert_q  =  \vert\vert t(u_2)u_2\vert\vert_q$, this implies that $t(u_1) = t(u_2)$ and $u_1 = u_2$. Therefore we have that $\varphi$ is one to one function. Let $w\in \Sigma_g$, $w\not= v$, set 
 $u = {w-v \over \vert\vert w - v\vert\vert_q}$, we have $t(u) =  \vert\vert w - v\vert\vert_q$ and   $\varphi({w-v \over \vert\vert w - v\vert\vert_q}) = w$. Thus, $\varphi$ is a bijection.
 
\noindent Suppose $\vert\vert v\vert\vert_q \geq 1$, let $\zeta\in C^{\infty}_c(\Omega)$ is such that 
$\vert\vert v- \zeta v\vert\vert_q < 1$. Observe that $ v- \zeta v = g$ on boundary. The same argument as above gives $t>0$ such that 
$\vert\vert v- t\zeta v\vert\vert_q =1$.
\cqfd
\section{Statement of further results}
 
\subsection{The sign of the Euler-Lagrange }

Let  $u$  be a minimizer for the problem (\ref{Pinf}), then, it  satisfies the following Euler-Lagrange equation
\ba \left\{
\begin{array}{lll} -\textsl{$\div$}(p(x)\nabla u)={ \Lambda u}^{q-{1}} &\textrm{in $\Omega$,}\\
\hspace{23mm}u > 0 &\textrm{in $\Omega$,}\\
\hspace{23mm}u=g  &\textrm{on $\partial \Omega$,}\\
\hspace{17mm}\vert \vert u \vert\vert_q =1,
\end{array}
\right. \label{Pequation} \ea 
where $ \Lambda \in \R$ is  the Lagrange multiplier associated to the problem (\ref{Pinf}), let $v$ be defined by (\ref{vvv}). The sign of $\Lambda$  is given by the  following:

\begin{theorem}\label{th2}
The sign of $\Lambda$  is as the following: 
If $\vert \vert v \vert\vert_q < 1$ then $\Lambda > 0$,
 if $\vert \vert v \vert\vert_q > 1$ then  $\Lambda < 0$
and if $\vert \vert v \vert\vert_q = 1$ then   $\Lambda = 0 $.
\end{theorem}

\subsection{Presence of a linear perturbation}
Over the course of the proof of Theorem \ref{th1}, one also reaps the following compactness result. 

\begin{theorem}\label{th3}
 We assume that $p$, $g$ and $v$ satisfy the same conditions as in Theorem \ref{th1}. Assume that $\vert\vert v\vert \vert_q < 1$. Let us denote by $\lambda_1$  the first eigenvalue of the operator $-\div(p \nabla.)$ with  homogeneous Dirichlet boundary condition. Then  for $\lambda <  \lambda_1$ we have the infimum in $S_{\lambda}(p,g)$ is achieved in the  following cases:
 \begin{enumerate} 
\item
$\lambda > 0$, $\alpha > 2$ and $n\geq 3$,
\item
$\lambda > 0$, $\alpha \leq 2$ and $n\in [3,2\alpha  +2[$.
\item
$\lambda < 0$, $n =3$ or $4$  with $\alpha > 1$ and $n =5$ with $\alpha > {3\over 2}$.
 \end{enumerate}
\end{theorem}
 
In the  presence of a linear perturbation, we will highlight a competition between three quantities, the dimension $n$,  the exponent $\alpha$ in (\ref{condition p}) and  the term of the linear perturbation. 
As we will see and as in Theorem \ref{th1}  the behavior of $p$ near its minimum plays an important role. The exponent $\alpha = 2$ is critical  in the case  $\lambda \not= 0$.

\section{Proof of Theorem 1.1.}

Let us start by proving the first part  of Theorem \ref{th1}. Suppose that $\vert\vert v\vert\vert_q < 1$. Since  the function  $u$ is a solution of $S_0(p,-g)$ if and only if $-u$ is solution of $S_0(p,g)$,
it  suffices to consider the case $g \geq 0$. 

\bigskip

Let $(u_j)$ be a minimizing sequence for $S_{0}(p,g)$, that is,

$$ \int_{\Omega} {p(x)}\vert\nabla{u_j(x)}\vert^{2}dx =  S_{0}(p, g) + o(1) $$
and

$$\vert\vert u_j \vert\vert_q = 1, \quad\quad u_j =  g \quad \hbox{in}\quad \partial\Omega.$$

\noindent Since $g\geq 0$, we may always assume  that $u_j \geq 0$, indeed, $(|u_j|)$ is also a minimizing sequence.  Since $(u_j)$ is bounded in $H^1$ we may extract a subsequence still  denoted by $(u_j)$ such that
$(u_j)$ converges weakly in $H^1$ to a function  $u\geq 0$ $a.e.$, $(u_j)$ converges  strongly to $u$ in $L^2 (\Omega)$, and $(u_j)$ converges to $u$  $a.e.$ on $\Omega$ with  $u =g$ on $\partial \Omega$.

Using a standard lower semicontinuity  argument,  we infer that  $\vert\vert u \vert\vert_q \leq 1$. To show that our infimum is achieved it suffices to prove that $\vert\vert u \vert\vert_q = 1$. Arguing by contradiction, let us assume that 
$$\vert\vert u \vert\vert_q <1.$$ 

We will prove that this is not possible with the assistance of several lemmas. We start by giving the first-order term of the energy $\int_{\Omega} {p(x)}\vert\nabla{u(x)}\vert^{2}dx$, next, we show that $u$ satisfies some kind Euler-Lagrange equation and then it is smooth. Finally, we compute the second-order term and highlight a contradiction.

\subsection{The first-order term}

\begin{lemma}\label{lemma1} 
For every $w\in H^1_g(\Omega)$ such that $\vert\vert w \vert\vert_q <1 $, we have

\begin{eqnarray}{\label{lemma1-1}}
S_{0}(p, g)  - \int_{\Omega} {p(x)}\vert\nabla{w(x)}\vert^{2}dx   \leq  p_0 S \left(1 - \int_\Omega \vert w\vert^q \right)^{2\over q} , 
\end{eqnarray}
For the weak limit $u$, we have equality: 
\begin{eqnarray}{\label{lemma1-2}}
S_{0}(p, g) -  \int_{\Omega} {p(x)}\vert\nabla{u(x)}\vert^{2}dx  =  p_0S\left(1 - \int_\Omega \vert u\vert^q \right)^{2\over q}.
 \end{eqnarray}
\end{lemma}

\Proof
Let  $w \in H^1_g(\Omega)$ such that  $\vert\vert w\vert\vert_q <1$. Therefore we can find a constant $c_{\varepsilon, a}  > 0$ such that 

\begin{eqnarray*}
\vert\vert w + c_{\varepsilon, a} u_{\varepsilon, a} \vert\vert_q  =1.
\end{eqnarray*}

\noindent Using Brezis-Lieb Lemma (see \cite{4}),  we obtain

\begin{eqnarray}{\label{cepsilon}}
c_{\varepsilon, a}^q = {1\over K_2 } \left(1 - \int_\Omega \vert w \vert^q \right) + o(1) 
\end{eqnarray}

\noindent where $K_2$ is defined in (\ref{K2}).
Careful expansion as $\varepsilon \rightarrow 0$ shows that  (see \cite{13}), for $n\geq 4$ 
\begin{equation}\label{uaepsilon}
\begin{array}{ll}
\hspace{-18mm}
\displaystyle\int_{\Omega}p(x)\vert{\nabla u_{{a},\varepsilon}(x)}\vert^{2}dx\leq\\[\bigskipamount] \left\{\begin{array}{llll}
\hspace{-1mm}p_{_{0}}K_{1}+O(\varepsilon^{n-2})&\hspace{-1.8mm}\textrm{if
$\left\{\begin{array}{lll}n\ge4\textrm{\quad and}\\ n-2< \alpha ,\end{array}\right.$}\\[\medskipamount]
\hspace{-1mm}p_{_{0}}K_{1}+A_{_{1}}\varepsilon^{\alpha}+o(\varepsilon^{\alpha})&\hspace{-1.8mm}\textrm{if $\left\{\begin{array}{lll}n\ge{4}\textrm{\quad and}\\n-2>\alpha,\end{array}\right.$}\\[\medskipamount]
\hspace{-1mm}p_{_{0}}K_{1}+\displaystyle  A_{_{2}}\ \varepsilon^{n-2} |\log\varepsilon  |+o(\varepsilon^{n-2}|\log\varepsilon|)&\hspace{-3.5mm}\textrm{
if $\left\{\begin{array}{ll}n \geq 4\textrm{\quad and}\\\ \alpha =n-2,
\end{array}\right.$}\\[\medskipamount]
\end{array}
\right. 
\end{array}
\end{equation}
\noindent with 
$$K_{1}=(n-2)^{2}\int_{\R^{n}}\frac{\vert y\vert^{2}}{(1+\vert
y\vert^{2})^{n}}dy$$
and  where $A_{1 }$,  $A_{2}$ and $A_{3}$ are positive constants depending only on $n$, $\gamma$ and $\alpha$, and for $n= 3$ and for $\alpha > 1$ we have as $\varepsilon \rightarrow 0$,
\begin{eqnarray*}
\int p(x)|\nabla u_{a,\varepsilon}(x)|^{2}dx = p_0 K_{1} + [\omega_{_{3}}\int_{0}^{R}(p_{0} + &\gamma r^{\alpha })|\psi'(r)|^{2}dr+\\
\omega_{_{3}}
k\alpha \int_{0}^{R}|\psi |^{2}r^{\alpha -2}dr] \varepsilon + o(\varepsilon). 
\end{eqnarray*}
where $\psi$ is defined as in (\ref{psi}).
\noindent Therefore for  $n= 3$ and  $\alpha > 1$ we obtain
\begin{equation}\label{n=3}
\int p(x)|\nabla u_{a,\varepsilon}(x)|^{2}dx =  p_0 K_{1} + A_{4}\varepsilon + o(\varepsilon). 
\end{equation}
where $A_{4}$ is a positive constant.

Remark that regardless of dimension $n$ as long as $n\geq 3$ and for $\alpha > 1$ we have 
\begin{eqnarray}\label{inferieur}
\int_{\Omega}p(x)\vert{\nabla  u_{{a},\varepsilon}(x)}\vert^{2}dx \leq  p_0K_1 +o(1).
\end{eqnarray}

\noindent Using $ w_\varepsilon = w + c_{\varepsilon, a} u_{\varepsilon, a}$ as testing function in  $S_{0}(p,g)$ we obtain
 \begin{eqnarray*}
S_{0}(p, g) \leq   \int_{\Omega} {p(x)}\vert\nabla{w(x)}\vert^{2}dx  +  c_{\varepsilon, a}^2
\int_{\Omega}p(x)\vert{\nabla  u_{{a},\varepsilon}(x)}\vert^{2}  + o(1) \end{eqnarray*}

\noindent Using (\ref{cepsilon}), the fact that ${K_1\over K_2^{2\over q}} = S$ and taking into account  (\ref{inferieur}) we get the first assertion of the Lemma \ref{lemma1}.\bigskip

\noindent For the second part, thanks to (\ref{lemma1-1}), it suffices to prove one inequality for $u$.
\begin{eqnarray}{\label{lemma=}}
S(p, g) -  \int_{\Omega} {p(x)}\vert\nabla{u(x)}\vert^{2}dx  \geq p_0S
\left(1 - \int_\Omega \vert u\vert^q \right)^{2\over q}.
 \end{eqnarray}

\noindent Set $v_j  = u_j -u$ so that  $v_j = 0 $ in $\partial \Omega$ and $(v_j)$ converges weakly to $0$ in $H^1_0$ and $a.e.$
We have by Sobolev inequality

\begin{eqnarray} \label{sobolev}
\int_\Omega p(x) \vert \nabla v_j\vert^2 \geq p_0 S\vert \vert v_j \vert \vert_q^2.
\end{eqnarray}

\noindent On the other hand,  we have (see \cite{4})

\begin{eqnarray}\label{BL}
1 = \int_\Omega \vert v_j\vert^q + \int_\Omega \vert u\vert^q + o(1).
\end{eqnarray}

\noindent Since $(u_j)$ is a minimizing sequence we have

\begin{eqnarray}\label{vj}
S_{0} (p,g) = \int_\Omega p(x) \vert\nabla  v_j\vert^2 + \int_\Omega p(x) \vert  \nabla u\vert^2 + o(1),
\end{eqnarray}

\noindent hence, combining (\ref{sobolev}), (\ref{BL}) and (\ref{vj}) we obtain the desired conclusion.\cqfd

We will now use the fact that $g$ is not identically zero. A consequence of the above lemma is the following:

\begin{lemma}\label{lemma2} 
The function $u$ satisfies
 \begin{eqnarray}\left\{\begin{array}{lll} \label{eulerequation}
 -\div(p\nabla u)& = p_0 S \left(1 - \int_\Omega \vert u\vert^q \right)^{2-q\over q}  \vert u\vert^{q-2}u &\textrm {in $\Omega$}\\[\medskipamount]
&u = g &\textrm{on $\partial \Omega$}
\end{array}
\right.
 \end{eqnarray}
Moreover, $u $ is smooth,  $u\in L^{\infty}(\Omega)$ and $u >0$  in $\Omega$.
\end{lemma}

\Proof
Applying (\ref{lemma1-1})   to  $w = u  + t\varphi$, $\varphi \in C^{\infty}_0(\Omega)$ and $\vert t \vert$ small enough, we have

\begin{eqnarray*}
S_{0} (p,g) \leq \int_\Omega p(x) \vert\nabla  u\vert^2 - 2t\int_\Omega p(x) \nabla u \nabla \varphi  +  o(t) + & \\
p_0 S \left(1 - \int_\Omega\vert u\vert^q - q t \int_\Omega\vert  u\vert^{q- 2}u \varphi + o(t)\right)^{{2 \over q}},
\end{eqnarray*}
thus
\begin{eqnarray*}
S_{0} (p,g) \leq \int_\Omega p(x) \vert\nabla  u\vert^2 - 2t\int_\Omega p(x) \nabla u \nabla \varphi  + &\\
p_0 S \left(1 - \int_\Omega\vert u\vert^q \right)^{{2 \over q}} \left( 1 - 2 t {\int_\Omega\vert  u\vert^{q- 2}u \varphi \over
1 - \int_\Omega\vert u\vert^q} + o(t)\right).
\end{eqnarray*}

\noindent Hence, by using  (\ref{lemma=}) we obtain for every $\varphi \in C^{\infty}_0(\Omega)$
\begin{eqnarray}\label{rien}
 - \int_\Omega p(x) \nabla u \nabla \varphi  - \left(1 - \int_\Omega\vert u\vert^q \right)^{2 -q \over q} \int_\Omega\vert  u\vert^{q- 2}u \varphi = 0.
\end{eqnarray}
 
\noindent Since $u=g$ on $\partial\Omega$ we obtain (\ref{eulerequation}).

For proving the regularity of $u$, it suffices, in view of the standard elliptic regularity theory  to show that $u$ is in $L^t(\Omega)$  for all $t<  \infty$. To see this, we shall apply Lemma  A1 of  \cite{5}, then, $u$ is 
as smooth as the regularity of $p$ and $g$ permits.

\noindent By using the strong  maximum  principle, and the fact that  $g \geq 0$, $g\not\equiv 0$ we get 
\begin{equation} \label{u>0}    
u >0  \quad \hbox{ in }\quad  \Omega.
\end{equation}
\cqfd

\subsection{The second-order term }
Now, we need   a refined version of (\ref{lemma1-1}). Similarly as in the proof of (\ref{lemma1-1}), let  $c_{\epsilon,a}$ be defined by
$1=\int_\Omega |u+c_{\epsilon,a}u_{\epsilon,a}|^q$. 
We can write  

\begin{eqnarray}\label{c0}
c_{\varepsilon, a}= c_0(1 - \delta(\varepsilon))
\end{eqnarray}

\noindent with

\begin{eqnarray}\label{c1}
c_0^q =  {1\over K_2 } \left(1 - \int_\Omega \vert u\vert^q \right) \quad\hbox{and}\quad \lim_{\varepsilon \rightarrow 0} \delta(\varepsilon) = 0.
\end{eqnarray}

\begin{lemma}\label{lemma2delta} 
We have
\begin{eqnarray} \label{delta}
\delta (\varepsilon) K_2 c_0^q &\geq& p_0 \varepsilon^{n-2 \over 2}
\left( c_0\int_\Omega  u^{q -1}{\psi \over \vert x - a \vert^{n-2}}c_0^{q}(q-1)D u(a) \right)\\ 
&+& 
\lefteqn{ {q-1\over 2}c_{0}^q K_{2} \delta^2(\varepsilon) + o(\delta^2(\varepsilon))+o(\varepsilon^{n-2\over 2}).}\nonumber 
 \end{eqnarray}
 where $D$ is a positive constant.
\end{lemma}
\Proof

First case $q\geq 3$. We need the following inequality,  for all $a \geq 0$ and $b\geq 0$ we have
\begin{eqnarray}\label{identity}
( a + b)^q \geq a^q + qa^{q-1}b + qab^{q-1} +b^q
\end{eqnarray}
which follows from 
$${ t^q + qt^{q-1} + qt +1 \over (1 + t)^q}\leq 1$$ 
for $t$ such that  $t ={b\over a}$ if $a \not= 0$.

\noindent Using  (\ref{identity}) and  the fact that $u >0$  we get  

\begin{eqnarray*}\label{inegalite} \displaystyle
1 &=& \int_\Omega \vert u + c_{\varepsilon, a} u_{\varepsilon, a} \vert^q\\
& \geq &
\int_\Omega  u ^q +
q c_{\varepsilon, a}^{q -1} \int_\Omega   u u_{\varepsilon, a}^{q -1} +
q c_{\varepsilon, a}\int_\Omega   u^{q -1} u_{\varepsilon, a} +
c_{\varepsilon, a}^q \int_\Omega     u_{\varepsilon, a}^q. 
\end{eqnarray*}

\noindent and thus

\begin{eqnarray}\label{inegalite} 
1 &\geq& \int_\Omega  u ^q +q c_{\varepsilon, a}^{q -1}
 \int_\Omega   u u_{\varepsilon, a}^{q -1} +
q c_{0} (1 - \delta(\varepsilon))\int_\Omega   u^{q -1} u_{\varepsilon, a} 
\\&+&\lefteqn{q c_{0}^q \left( 1 -q\delta(\varepsilon) +{ q(q-1)\over 2} \delta^2(\varepsilon)+ o( \delta^2(\varepsilon)\right)  \int_\Omega     u_{\varepsilon, a}^q.}\nonumber 
\end{eqnarray}

\noindent On the other hand we have

\begin{eqnarray}\label{dirac}\displaystyle
\int_\Omega   u u_{\varepsilon, a}^{q -1} = \varepsilon^{n - 2 \over 2} D u(a) + o(\varepsilon^{n - 2 \over 2})
\end{eqnarray}

\noindent where $D$ is a positive constant, and

\begin{eqnarray}\label{dirac1}
 \int_\Omega   u^{q -1}  u_{\varepsilon, a}   = \varepsilon^{n - 2 \over 2} \int_\Omega u^{q -1}  {\psi \over \vert x- a \vert^{n- 2}}
+ o(\varepsilon^{n - 2 \over 2}).
\end{eqnarray}
\noindent Combining  (\ref{K2}), (\ref{inegalite}), (\ref{dirac})  and (\ref{dirac1}) we obtain  (\ref{delta}).
\bigskip

\noindent Second case $2 < q <3$. In what follows $C$ denote a positive constant independent of $ \varepsilon$.  
The keys are  the two   following inequalities,  we have for all $a \geq 0$ and $b\geq 0$ 

\begin{eqnarray}\label{identity2}
\vert ( a + b)^q -( a^q + qa^{q-1}b + qab^{q-1} +b^q)\vert \leq Ca^{q-1}b \quad \hbox{if}\,  a \leq b
\end{eqnarray}
\noindent and
\begin{eqnarray}\label{identity3}
\vert ( a + b)^q -( a^q + qa^{q-1}b + qab^{q-1} +b^q)  \vert \leq Cab ^{q-1}\quad\hbox{if} \, a \geq b
\end{eqnarray}
which follows respectively  from 
\begin{eqnarray}
{\vert (1 + t)^q -( t^q + qt^{q-1} + qt  +1) \vert \over t} \leq C
\end{eqnarray} 
for $t\geq 1$ and 
\begin{eqnarray} {\vert (1 + t)^q -( t^q + qt^{q-1} + qt  +1) \vert \over t^{q-1}} \leq C
\end{eqnarray}
for $t\leq 1$  for $t$ such that   $t ={b\over a}$ if $a \not= 0$.
\bigskip

\noindent Using  (\ref{identity2}) and (\ref{identity3}) we get  
\begin{eqnarray} \displaystyle
1 &=& \int_\Omega \vert u + c_{\varepsilon, a} u_{\varepsilon, a} \vert^q  \label{egalite}
\\&=&
\int_\Omega  u ^q + q c_{\varepsilon, a}^{q -1} \int_\Omega   u u_{\varepsilon, a}^{q -1} +
q c_{\varepsilon, a}\int_\Omega   u^{q -1} u_{\varepsilon, a} +
c_{\varepsilon, a}^q \int_\Omega     u_{\varepsilon, a}^q 
 \nonumber \\&+& 
\lefteqn{R^{(1)}_\varepsilon + R^{(2)}_\varepsilon.}\nonumber
\end{eqnarray}

\noindent where 
\begin{eqnarray*}\label{reste1} 
R^{(1)}_\varepsilon  \leq C\int_{\{x,  u\geq c_{\varepsilon, a}\psi U_{a, \varepsilon} \}}u \vert \psi U_{a, \varepsilon}\vert^{q-1}
\end{eqnarray*}
\noindent and

\begin{eqnarray*}
R^{(2)}_\varepsilon  \leq C\int_{\{x, u < c_{\varepsilon, a}\psi U_{a, \varepsilon}\} } u^{q-1} \psi U_{a, \varepsilon}.
\end{eqnarray*}

\noindent We claim that the remainders terms $R^{(1)}_\varepsilon $ and $R^{(2)}_\varepsilon $ verify
\begin{eqnarray}\label{restes} 
R^{(1)}_\varepsilon  = o(\varepsilon^{n-2\over 2})\quad \hbox{and} \quad R^{(2)}_\varepsilon =  o(\varepsilon^{n-2\over 2}).
\end{eqnarray}

\noindent Let us justify the first assertion in (\ref{restes}).
In the set $\Omega\setminus  B(a,r)$ we have $U_{a, \varepsilon}^{q-1} \leq 
C \varepsilon^{n + 2\over 2}$ and in the set 
$B(a, r) \cap \{x, u\geq c_{\varepsilon, a}\psi U_{a, \varepsilon} \}$ we have  
 $U_{a, \varepsilon} \leq C$ and then necessarily $\vert x - a \vert \geq C\varepsilon^{1\over 2}$, therefore

\begin{equation}\label{reste1-1} 
R^{(1)}_\varepsilon  \leq C\int_{\{x,  C\varepsilon^{1\over 2} < \vert x - a \vert \leq r \} }
\left( {\varepsilon \over \varepsilon^2 + \vert x - a \vert^2}\right)^{n + 2\over 2} dx 
=  o(\varepsilon^{n-2\over 2}).
\end{equation}
\bigskip
\noindent Let us verify that   $R^{(2)}_\varepsilon =  o(\varepsilon^{n-2\over 2})$. In the set $ A_{a, \varepsilon} = \{x, u <  c_{\varepsilon, a}\psi 
 U_{a, \varepsilon} \}$ we have $\psi > 0$ and consequently, since $u$ is smooth, there exists $\delta > 0$ such that $u > \delta $ in $A_{a, \varepsilon} $ thus $U_{a, \varepsilon} \geq C $. This implies that
$\vert x - a \vert \leq C \varepsilon^{1\over 2}$. We have 

\begin{equation}\label{reste1-2} 
R^{(2)}_\varepsilon  \leq C\int_{\{x, \vert x - a \vert \leq C \varepsilon^{1\over 2} \} }
\left( {\varepsilon \over \varepsilon^2 + \vert x - a \vert^2}\right)^{n - 2\over 2} dx  =  o(\varepsilon^{n-2\over 2}).
\end{equation}
Combining    (\ref{egalite}), (\ref{dirac}), (\ref{dirac1}),  and  (\ref{restes}) we obtain  that $ \delta(\varepsilon) = O(\varepsilon^{n-2\over 2})$ and  (\ref{delta}).
\cqfd
\noindent We are able to prove now:
\begin{lemma}\label{lemma3} 
If $n\geq 3$ and $\alpha >1$ then we have for every $ 3 \leq n < 2\alpha + 2$ we have 
\begin{eqnarray}\label{contradiction1}
 S_{0}(p, g)-  \int_{\Omega} {p(x)}\vert\nabla{u(x)}\vert^{2}dx   < p_0S\left(1 - \int_\Omega \vert u\vert^q \right)^{2\over q}.  
 \end{eqnarray}
\end{lemma}
Let us postpone the proof of Lemma \ref{lemma3} and complete the first part of the  proof  of Theorem \ref{th1}. Combining (\ref{contradiction1}) and (\ref{lemma1-2})  this  leads to a contradiction and then we obtain that $\vert\vert u \vert\vert_q =1$ and  therefore the infimum $S_{0}(p, g)$ is achieved.
\Proof of the first part of Theorem \ref{th1}.
Let us chose  $ w_\varepsilon = u + c_{\varepsilon, a} u_{\varepsilon, a}$  as testing function in  $S_{0}(p,g)$, we obtain   
\begin{eqnarray}\label{inf} 
S_{0}(p, g) \leq \int_\Omega p\vert \nabla (u + c_{\varepsilon, a} u_{a,\varepsilon}) \vert^2.
\end{eqnarray}
\noindent By  (\ref{inf}) and  (\ref{c0}) it is easy to  see
\begin{eqnarray*}
S_{0}(p, g) &\leq& 
\int_\Omega p\vert \nabla u\vert^2 - 2 c_0\varepsilon^{n-2\over 2} \int_\Omega 
(\div(p\nabla u)  {\psi \over \vert x  - a \vert^{n-2}} \\&+&
c_0^2(1 - 2\delta( \varepsilon)  + \delta^2(\varepsilon) )  \int_{\Omega} p(x) |\nabla u_{a,\varepsilon}|^{2}dx  + o(\varepsilon^{n-2\over 2}). 
\end{eqnarray*}

\noindent Now, using (\ref{delta}) and the fact that  $\delta(\varepsilon) = o(1)$ we infer
\begin{eqnarray} \label{devel-delta1}
&S_{0}(p, g)\leq
 \int_\Omega p\vert \nabla u\vert^2  + p_{0} K_1 c_0^2 -  2 c_0 \varepsilon^{n-2\over 2}  \int_\Omega  \div(p\nabla u)  {\psi \over \vert x  - a \vert^{n-2}}\\
 \lefteqn{ - 2 c_0^{2}   \left[ {\varepsilon^{n-2 \over 2} \over K_2 c_0^q}  \left( c_0\int_\Omega  u^{q -1}{\psi \over \vert x - a \vert^{n-2}} + c_0^{q-1}D u(a) \right) +
  {q-1\over 2}\delta^2(\varepsilon) + o(\delta^2(\varepsilon))  \right]} \nonumber  \\
 \lefteqn{ \int_{\Omega} p(x) \vert \nabla u_{a,\varepsilon} \vert^{2}dx
 +c_0^{2}\delta^2(\varepsilon)   \int_{\Omega} p(x) \vert \nabla u_{a,\varepsilon} \vert^{2}dx  + o(\varepsilon^{n-2\over 2}).}\nonumber
\end{eqnarray}

\noindent Since $  \int_{\Omega} p(x) \vert \nabla u_{a,\varepsilon} \vert^{2}dx = K_2 + o(1)$ we obtain

\begin{eqnarray}\  \label{devel-delta2}
&S_{0}(p,g) \leq \int_\Omega p\vert \nabla u\vert^2  +c_0^2 \int_{\Omega} p(x) |\nabla u_{a,\varepsilon}|^{2}dx   -(q-2)c_0^2 \delta^2(\varepsilon) + o(\delta^2(\varepsilon)) -\nonumber 
\\ \lefteqn{2 c_0\left[ \int_\Omega  \div(p\nabla u)  {\psi \over \vert x  - a \vert^{n-2}}
 + \left( {c_{0}^{2-q} \over K_2 }  \int_\Omega  u^{q -1}{\psi \over \vert x - a \vert^{n-2}}  
 + { D\over K_{2}} u(a)\right)   (K_2 + o(1))\right]\varepsilon^{n-2\over 2}}\nonumber \\ 
\lefteqn{ +  o(\varepsilon^{n-2\over 2}).}\nonumber
\end{eqnarray}
\noindent This leads to
\begin{eqnarray*}  
S_{0}(p, g) &\leq& \int_\Omega p\vert \nabla u\vert^2  +c_0^2 \int_{\Omega} p(x) |\nabla u_{a,\varepsilon}|^{2}dx -
(q-2)c_0^2K_{2}\delta^2(\varepsilon) + o(\delta^2(\varepsilon))\\ &-&
\lefteqn{2c_{0}{ DK_{1}\over K_{2}} u(a)  \varepsilon^{n-2\over 2}  + o(\varepsilon^{n-2\over 2}).}\nonumber
\end{eqnarray*}

\noindent We know that  $\delta^2(\varepsilon)= o(1)$ thus
\begin{eqnarray} \label{devel-avec-delta}
S_{0}(p, g) &\leq& \int_\Omega p\vert \nabla u\vert^2  +c_0^2 \int_{\Omega} p(x) |\nabla u_{a,\varepsilon}|^{2}dx\\ &-&
2c_{0}{ DK_{1}\over K_{2}} u(a)  \varepsilon^{n-2\over 2}  + o(\varepsilon^{n-2\over 2})\nonumber.
\end{eqnarray}

\noindent We  are now able  to give a precise asymptotic behavior of the RHS of (\ref{inf}). This will possible thanks to the fact that $u(a)\neq 0$, namely $u(a) >  0$. One needs to distinguish between dimensions and the parameter $\alpha$. 
Four cases  follow from  (\ref{uaepsilon}) and  (\ref{devel-avec-delta}):
\begin{itemize}
\item
The case when $n\geq 4$ and $n <  \alpha + 2$. We have
\end{itemize}

\begin{eqnarray*}\  \label{devel-delta2}
S_{0}(p, g) \leq \int_\Omega p\vert \nabla u\vert^2  + c_0^2 \left(  p_{_{0}}K_{1} + o(\varepsilon^{n - 2 })\right)   - 2c_{0}  { DK_{1}\over K_{2}} u(a)  \varepsilon^{n-2\over 2}  + o(\varepsilon^{n-2\over 2}).
\end{eqnarray*}
\noindent Consequently,  we have

\begin{equation}\label{1}
S_{0}(p, g)  \leq \int_\Omega p\vert \nabla u\vert^2 +  p_0 c_0^2 K_{1} - 2c_{0}  { DK_{1}\over K_{2}} u(a) \varepsilon^{n-2\over 2}  + 
o(\varepsilon^{n-2\over 2}).
\end{equation}
\begin{itemize}
\item
The case when $n\geq 4$ and $ n  > \alpha + 2$. We have
\end{itemize}

\begin{eqnarray*}\  \label{devel-delta2}\displaystyle
S_{0}(p, g) &\leq& \int_\Omega p\vert \nabla u\vert^2  +c_0^2 \left(  p_{_{0}}K_{1} + A_{_{2}} \varepsilon^{\alpha} + o(\varepsilon^{\alpha})\right) \\  &-&
 2c_{0}  { DK_{1}\over K_{2}}  u(a)  \varepsilon^{n-2\over 2}  + o(\varepsilon^{n-2\over 2}).
\end{eqnarray*}

\noindent Therefore,  we have
\begin{eqnarray*}\displaystyle
S_{0}(p, g)  \leq \int_\Omega p\vert \nabla u\vert^2  +  p_0 c_0^2 K_{1} - 2c_{0}  { DK_{1}\over K_{2}} u(a) \varepsilon^{n-2\over 2} + A_{2}c_0^2 \varepsilon^\alpha +   o(\varepsilon^\alpha) + 
o(\varepsilon^{n-2\over 2}).
\end{eqnarray*}
Hence, if $n < 2\alpha + 2$ then
 \begin{eqnarray}\label{2}
S_{0}(p, g)  \leq& \int_\Omega p\vert \nabla u\vert^2  - 2c_{0}  { DK_{1}\over K_{2}} u(a) \varepsilon^{n-2\over 2}  + 
o(\varepsilon^{n-2\over 2}).
\end{eqnarray}

\begin{itemize}
\item
The case when $n \geq  4$ and $\alpha  = n - 2$. We have
\end{itemize}

\begin{eqnarray*}\  \label{ddevel-avec-delta1}
S_{0}(p, g) &\leq& \int_\Omega p\vert \nabla u\vert^2  +c_0^2 \left(p_{_{0}}K_{1}+\displaystyle  A_{_{2}}\ \varepsilon^{n-2} |\log\varepsilon  |+o(\varepsilon^{n-2}|\log\varepsilon|)\right) \\
&-& 2c_{0} { DK_{1}\over K_{2}} u(a)  \varepsilon^{n-2\over 2}  + o(\varepsilon^{n-2\over 2}).
\end{eqnarray*}
\noindent thus we  get

\begin{equation}\label{3}
 S_{0}(p, g)  \leq \int_\Omega p\vert \nabla u\vert^2  + p_0 c_0^2K_1 - 2c_{0} { DK_{1}\over K_{2}}u(a) \varepsilon^{n-2\over 2} +  o(\varepsilon^{n-2\over 2}).
\end{equation}

\begin{itemize}
\item{4}
The case when $n = 3$ and $\alpha > 1$. We have
\end{itemize}

\begin{eqnarray*}\  \label{devel-delta2}
S_{0}(p, g) \leq \int_\Omega p\vert \nabla u\vert^2  + c_0^2  [ p_0 K_1 + A_{4}\varepsilon + o(\varepsilon)] -
2c_{0}  { DK_{1}\over K_{2}} u(a)  \varepsilon^{1\over 2}  + o(\varepsilon^{1\over 2}).
\end{eqnarray*}
\noindent Hence  we have

\begin{equation}\label{4}
S_{0}(p, g) \leq \int_\Omega p\vert \nabla u\vert^2  +  p_0 c_0^2 K_{1} -  2c_{0}  
{ DK_{1}\over K_{2}}  u(a) \varepsilon^{1\over 2} +  o(\varepsilon^{1\over 2}).
\end{equation}

\noindent Now, thanks to (\ref{1}), (\ref{2}), (\ref{3}), (\ref{4})  and the fact that $u(a) >0$  we obtain the estimates in Lemma \ref{lemma3}.
\cqfd

\subsection{The case $\vert\vert v\vert\vert_q \geq 1$.}
For the proof of  the second part of Theorem \ref{th1}  we set
\begin{equation*}
\alpha := \inf_{{u}\in
{H}^{1}_{g}(\Omega),\|u\|_{q}=1}\int_{\Omega}{p(x)}\vert\nabla{u(x)}\vert^{2}dx
\end{equation*}
and 
\begin{equation*}
\beta := \inf_{{u}\in
{H}^{1}_{g}(\Omega),\|u\|_{q}\leq1}\int_{\Omega}{p(x)}\vert\nabla{u(x)}\vert^{2}dx.
\end{equation*}

\noindent Indeed using the convexity of the problem $\beta$, it is clear that  the infimum in $\beta$ is achieved by some function  
$w \in {H}^{1}_{g}(\Omega)$ satisfying $\vert\vert w \vert\vert_q \leq1$. 
Necessarily we have  equality.  Let us reason by contradiction, if we had $\vert\vert w \vert\vert_q < 1$, 
let $\zeta\in ^\infty_{c}(\Omega)$, for $t$ real and small such that we have $\vert\vert w + t\zeta \vert\vert_q < 1$, using $ w + t\zeta $ as test function in $\beta$ we obtain that  $w$ would  be the unique solution of the following Euler-Lagrange equation:
\begin{equation}
\left\{
\begin{array}{llllll}
-\div(p\nabla w) =&0 \textrm{\quad in $\Omega$},\\
w=g&  \textrm{\quad on $\partial\Omega$}.
\end{array}
\right.
\label{equationg}
\end{equation}

\noindent that is mean $w$ and $v$ coincide, this leads to a contradiction since $\vert\vert v \vert\vert_q \geq 1$. Therefore $\alpha$ is achieved.

\noindent Since $\vert\vert w \vert\vert_q =1$ we have   $\int_{\Omega}{p(x)}\vert\nabla{w(x)}\vert^{2}dx = \beta \leq \alpha \leq \int_{\Omega}{p(x)}\vert\nabla{w(x)}\vert^{2}dx$. Thus $\alpha = \beta$.
\cqfd
\section{ The sign of the Euler-Langange multiplier. Proof of Theorem \ref{th2}}

We follow an idea of \cite{11}. Let  $u$  be a minimizer for the problem (\ref{Pinf}) and  $v$ be defined by (\ref{vvv}), using the fact that problem (\ref {vvv}) has a unique solution which minimizes (\ref{Pv}), we remark that we have  $\vert \vert v \vert\vert_q \not= 1$ if and only if we have $\Lambda \not= 0$.

\noindent Using (\ref{Pequation}) and (\ref{vvv}) we obtain
\ba \left\{
\begin{array}{lll} -\textsl{$\div$}(p(x)\nabla (u - v)) = { \Lambda u}^{q-{1}} &\textrm{in $\Omega$,}\\
\hspace{23mm}u-v=0  &\textrm{on $\partial \Omega$.}
\end{array}\label{u-v}
\right. \label{eq(u-v)} \ea

\noindent First, suppose  that $\vert \vert v \vert\vert_q < 1$. Multiplying (\ref{eq(u-v)}) by  $u-v$ and integrating we obtain 
\begin{equation}
 \Lambda (\vert \vert u \vert\vert_q^q   -\int_{\Omega} \vert u\vert^{q-1}v)   = \int_{\Omega} p(x)\vert \nabla (u - v)\vert^{2}.  \label{green} 
\end{equation}
\noindent From  H\"older inequality and the fact that $\vert\vert u \vert\vert_q = 1$ we obtain 
\begin{equation}
\vert \vert u \vert\vert_q^q   -\int_{\Omega} \vert u\vert^{q-1}v  \geq  1 -  \vert \vert v \vert\vert_q >0
 \label{holder} .
\end{equation}

\noindent Putting together  (\ref{green})  and  (\ref{holder}) and using the fact that $u \not= v$ we see that $ \Lambda >0$.
\noindent Suppose now that $\vert \vert v \vert\vert_q >1$. For $t\in \R$, let us define the function $f$ by 
\begin{equation*}
f(t) = \int_{\Omega} \vert tu + (1-t)v\vert^q dx.
\end{equation*}
\noindent Note  that the function $f$ is smooth and convex since $f''(t) =  q(q-2) \int_\Omega   tu + (1-t)v\vert^{q-1}(u-v)^{2} \geq 0$ and we have 

\begin{equation} \label{f(0) f(1)}
f(0) = \vert \vert v \vert\vert_q^q > 1 \quad\hbox{ and  } \quad  f(1) = \vert \vert u\vert\vert_q^q = 1. 
\end{equation}

\noindent We may use the following:

\begin{lemma}\label{lemma4} 
For all $t\in [0,1[$ we have $f(t) > 1$. 
\end{lemma}

\Proof 
Arguing by contradiction, since $f$ is continuous, by  the intermediate value theorem there exists $t_{0} \in [0,1[$ such that $f(t_0) = 1$. 
Using $t_{ 0}u + (1-t_{0})v \in \Sigma_{g}$ as testing function in $S_{0}(,p,g)$ we have
\begin{equation} \label{testingfunction}
S_{0}(p,g) = \int_{\Omega} p(x)\vert \nabla  u\vert^{2}  \leq  \int_{\Omega} p(x)\vert \nabla ( t_{ 0}u + (1-t_{0})v\vert^{2}
\end{equation}

\noindent Multiplying (\ref{vvv})  by $u-v$ and integrating we obtain
\begin{equation}\label{multiplyv}
\int_\Omega  p\vert \nabla v \vert^{2} = \int_\Omega p\nabla u  \nabla v 
\end{equation}

\noindent Using (\ref{testingfunction}),  (\ref{multiplyv})  and the fact that $t_{0} < 1$ we obtain
\begin{equation*}
\int_\Omega p\vert \nabla u \vert^{2} \leq \int_\Omega p\vert \nabla v  \vert^2
\end{equation*}
\noindent Since $v$ is the unique solution of (\ref{vvv}) we obtain that $u = v$ which clearly contradicts (\ref{f(0) f(1)}). This complete the proof of Lemma \ref{lemma4}.
\cqfd

\noindent By the convexity of $f$ and Lemma \ref{lemma4} we deduce that  $f'(1) \leq 0$. But 
$f'(1)  = q\int_{\Omega} \vert u\vert^{q-1} (u-v)$ and then by (\ref{u-v}) we have 
$f'(1)  = {q\over \Lambda} \int_{\Omega} p(x)\vert \nabla (u - v)\vert^{2}$. We conclude that $\Lambda <0$.
\cqfd
\section{Existence of minimizer in the presence of a linear perturbation:
Proof of  Theorem \ref{th3}}
First, we claim that if problem (\ref{Pinflambda}) has a solution then $\lambda < \lambda_1$. Indeed, let $u$ be a solution of (\ref{Pinf}) and $v$ satisfying (\ref{vvv}), we have

\ba \left\{
\begin{array}{lll} -\textsl{$\div$}(p(x)\nabla ( u - v )= \Lambda(\lambda , u) {u}^{q-{1}}  + \lambda u&\textrm{in $\Omega$,}\\
\hspace{30mm}u > 0 &\textrm{in $\Omega$,}\\
\hspace{23mm}u - v =0  &\textrm{on $\partial \Omega$.}
\end{array}
\right. \label{Pequation v} \ea 
\noindent where $\Lambda(\lambda , u)$ is a Euler-Lagrange multiplier. Since $\vert\vert v\vert\vert_q < 1$, using section 5, we find that $\Lambda(\lambda , u) > 0$. Let $\varphi_1$  be the  eigenfunction  of the operator $-\div(p \nabla.)$ with  homogeneous Dirichlet boundary condition corresponding to $\lambda_1$. Multiplying (\ref{Pequation v}) by $\varphi_1$ and integrating we obtain 
\begin{eqnarray*}\displaystyle
-\int_\Omega \textsl{$\div$}(p(x)\nabla (u - v) )\varphi_1 &=&
 \lambda_1 \int_\Omega (u - v)\varphi_1 \\ &=& \Lambda(\lambda , u)\int_\Omega {u}^{q-{1}} \varphi_1  + \lambda \int_\Omega u \varphi_1.
\end{eqnarray*}
Then we get 
\begin{eqnarray*}\displaystyle
 (\lambda_1 -  \lambda)\int_\Omega (u - v)\varphi_1 \geq \lambda_1\int_\Omega v\varphi_1
\end{eqnarray*}

\noindent and thus $\lambda < \lambda_1$.
 
The proof of Theorem \ref{th1}  is similar to the one of Theorem \ref{th1} so that we briefly outline it. We need only to take into account  the linear perturbation term. We will then follow exactly all the steps in the proof of Theorem \ref{th1} untill (\ref{devel-avec-delta}), we just need to account the linear perturbation. We get
 \begin{eqnarray}\  \label{v lambda}
S_{\lambda}(p, g) \leq &\int_\Omega p\vert \nabla u\vert^2+c_0^2 \left( \int_{\Omega} p(x) |\nabla u_{a,\varepsilon}|^{2}dx 
-\lambda  \int_{\Omega} |u_{a,\varepsilon}|^{2}dx\right) \\
-\lefteqn{ 2c_{0}{ DK_{1}\over K_{2}} u(a)  \varepsilon^{n-2\over 2}  + o(\varepsilon^{n-2\over 2}).}\nonumber
\end{eqnarray}
From  \cite{6} we have
\ba {\parallel
u_{a,\varepsilon}\parallel^{2}_{2}}= \left\{
\begin{array}{l} K_{3}\varepsilon^2+O(\varepsilon^{n-2}) \quad
i f\quad n\geq{5},
\\[\medskipamount]
C_1 \varepsilon^2 | \log {\varepsilon }| + O(\varepsilon^2)\quad \ if \quad n={4},
\\[\medskipamount]
C_2\varepsilon + O(\varepsilon^2)\ \quad if\quad  n={3}
\end{array}
\right.\label{L2} \ea 
\noindent where  $C_1$, $C_2$ and $C_3$ are positive constants.  
Using (\ref{uaepsilon}), (\ref{v lambda}) and (\ref{L2}) and the fact that $ u(a) >0$ we conclude the proof of
 Theorem \ref{th3}.
\cqfd

\end{document}